\documentclass[12pt,a4paper]{article}

\usepackage{amssymb,amsmath,amsthm,xspace,epsfig, amsfonts,amscd}

\usepackage[left=1.5cm, right=1.5cm, top=1.5cm]{geometry}

\begin{document}

\newcommand{\N}{\mathbb{N}}
\newcommand{\R}{\mathbb{R}}
\newcommand{\Z}{\mathbb{Z}}
\newcommand{\Q}{\mathbb{Q}}
\newcommand{\C}{\mathbb{C}}
\newcommand{\PP}{\mathbb{P}}

\newcommand{\LL}{\Bbb L}
\newcommand{\OO}{\mathcal{O}}

\newcommand{\esp}{\vskip .3cm \noindent}
\mathchardef\flat="115B

\newcommand{\lev}{\text{\rm Lev}}

\def\ut#1{$\underline{\text{#1}}$}
\def\CC#1{${\cal C}^{#1}$}
\def\h#1{\hat #1}
\def\t#1{\tilde #1}
\def\wt#1{\widetilde{#1}}
\def\wh#1{\widehat{#1}}
\def\wb#1{\overline{#1}}

\def\restrict#1{\bigr|_{#1}}

\def\hu#1#2{\mathsf{U}_{fin}\bigl({#1},{#2}\bigr)}
\def\ch#1#2{\left(\begin{array}{c}#1 \\ #2 \end{array}\right)}

\newtheorem{lemma}{Lemma}[section]

\newtheorem{thm}[lemma]{Theorem}
\newtheorem*{thm*}{Theorem}
\newtheorem*{lemma?}{Lemma ??}

\newtheorem{defi}[lemma]{Definition}
\newtheorem{conj}[lemma]{Conjecture}
\newtheorem{cor}[lemma]{Corollary}
\newtheorem{prop}[lemma]{Proposition}
\newtheorem{prob}[lemma]{Problem}
\newtheorem{qu}[lemma]{Question}
\newtheorem{q}[lemma]{Questions}
\newtheorem*{rem}{Remark}
\newtheorem{examples}[lemma]{Examples}
\newtheorem{example}[lemma]{Example}

\title{(Non-)contractible road spaces of trees}
\date{\today}
\author{Mathieu Baillif}
\maketitle

\abstract{\footnotesize 
   We prove that the road space of an $\R$-special tree is contractible and that a locally metrizable space containing
   a copy of an uncountable $\omega_1$-compact subspace of a tree is not. 
   We also
   raise some questions about possible generalizations.
   }

\section{Introduction and definitions}

Despite the word `contractible' in our title, our methods and results belong to set theoretic topology and have very little
in common with homotopy theory.
By `space' is meant a topological Hausdorff space, in particular `regular' and `normal' imply
Hausdorff.
We assume some familiarity with set theoretical trees
and countable ordinals. Our terminology is standard, definitions are given at the end of this introduction if needed.
Any tree is endowed with the order topology (also called the interval topology).

Asking questions about the contractibility of set theoretic objects is not very common, so before stating our results let us explain why 
we think it is an interesting subject, particularly for manifolds and similar spaces. 
(The main reason is that we like to think about these problems and find it oddly satisfying to try to mentally squish things, 
but we fear it will not be seen as valuable.)

The objects we study in this note are set theoretic trees and their road spaces, 
which are obtained by joining consecutive points by a line segment with a topology that makes it locally 
embeddable in $\R^2$ (the details are given below).
These spaces are a good toy model for so-called Type I non-metrizable manifolds (also defined below).
Metrizable manifolds that are contractible happen to coincide with those having vanishing homotopy groups:
Milnor showed in \cite{MilnorCW} that a metrizable manifold has the homotopy type of a CW-complex and  
Whitehead Theorem applies
(see, e.g., \cite[p. 346]{Hatcher}).
If one releases the metric assumption, this is not true anymore as
there are simple non-metrizable non-contractible manifolds with vanishing homotopy groups (e.g. the {\em longray} described below).
Of course, any space with non-trivial homotopy groups is non-contractible.
It was at first not clear to us whether contractible non-metrizable manifolds do exist at all.
We then noticed that examples were found a long time ago by Calabi and Rosenlicht;
the Pr\"ufer surface and some of its variants are contractible.
(The original source is \cite{CalabiRosenlicht} and a more recent account is given in \cite{mesziguesBridges}, see also \cite{GabardWouuuh} for more
on Pr\"ufer surface and homotopy.) 
This particular manifold has the property that
some open Euclidean set has a non-metrizable closure, in a terminology dating to P. Nyikos \cite{Nyikos:1984} it is in the class of
{\em Type II} manifolds. 
Hence, in a loose sense, what makes it non-metrizable
lies just at the boundary of some perfectly nice Euclidean open contractible set. We happen to be able to 
`push' this non-metrizable stuff inside this Euclidean open set, everything at once,
yielding the contractibility.
The other class of non-metrizable manifolds are the aforementioned {\em Type I} manifolds \cite[Def. 2.10]{Nyikos:1984}.
A space $X$ is of Type I if and only if $X=\cup_{\alpha\in\omega_1}X_\alpha$, where $X_\alpha$ is open and $\wb{X_\alpha}$ Lindel\"of
for each $\alpha$, and $\wb{X_\alpha}\subset X_\beta$ whenever $\alpha<\beta<\omega_1$ (and of Type II otherwise).
In case of manifolds, each $X_\alpha$ is an open metrizable submanifold since Lindel\"ofness of $\wb{X_\alpha}$ implies its metrizability.
Loosely speaking again, these manifolds `grow slowly' 
instead of jumping at once into non-metrizability. 
Since the Pr\"ufer manifold is of Type II, we are 
interested in finding whether there are contractible non-metrizable Type I (slowly growing) manifolds.
As written above, road spaces of trees are good toy models. 
Our first result gives a contractible example in this class of spaces:

\begin{thm}\label{thm:main1}
   The road space $R_T$ of a rooted $\R$-special tree is contractible.
\end{thm}

We shall show in a subsequent paper that given an $\omega_1$-tree $T$, it is possible to define
a surface which is is homotopy equivalent to $R_T$ (and contains it), providing the example of a contractible Type I manifold we are seeking.

Looking at the problem from the other side, it seemed interesting as well
to see whether there are properties purely from general topology 
(which have nothing to do with homotopy/homology groups or such) 
which prevent a manifold (or a `locally nice' space) to be contractible.
Some results are already available, for instance
the following theorem was proved by S. Deo and D. Gauld (using ideas developed by the author in \cite{meszigueshom}):

\begin{thm}[{\cite[Thm 3.4]{DeoGauld:2007}}]\label{thmdeogauld}
Any locally metrizable space $X$ containing a copy of $\omega_1$ is non-contractible.
\end{thm}

This implies for instance that any tree with an uncountable branch has a non-contractible road space.
An older and simpler result (see for instance \cite[Prop. 1.22]{GauldBook}) is that the longray defined as
inserting a line segment between consecutive ordinals in $\omega_1$ 
(more formally, $\omega_1\times[0,1)$ with the lexicographic order topology) is 
a non-contractible manifold with vanishing homotopy groups. 

We generalize Deo and Gauld's result in two ways. 
The first is to weaken the assumption that the subspace is a copy of $\omega_1$.

\begin{thm}\label{thm:main2}
   Let $T$ be a tree, $S\subset T$ be uncountable and $\omega_1$-compact in the subspace topology, 
   and let $X$ be a locally metrizable space.
   Then there is no continuous $h\colon S\times[0,1]\to X$ 
   such that $h(\cdot,1)\colon S\to X$ has uncountable image and $h(\cdot,0)$ is constant.
   In particular, if $X$ contains a copy of $S$ then $X$ is not contractible.
\end{thm}

Recall that a space is {\em $\omega_1$-compact} if and only if any closed discrete subset is at most countable.
Notice that an uncountable $\omega_1$-compact subset of a tree cannot be metrizable, since $\omega_1$-compactness
is equivalent to Lindel\"ofness in metrizable spaces (see, e.g., \cite[Thm 4.1.15]{Engelking}), 
and a Lindel\"of subset of a tree is countable.
We note that we could add the assumption that $T$ has height $\omega_1$ without loss of generality,
because if $S$ contains an uncountable branch $E$, then $E$ is $\omega_1$-compact in the subspace topology and thus homeomorphic to a 
stationary subset of $\omega_1$, and we may apply Theorem \ref{thm:main3} below.
Theorem \ref{thm:main2} gives us another motivation for these homotopy questions: to
find where is the line among the set theoretic or topological properties of 
trees between contractibility and non-contractibility of their road spaces
(see Question \ref{q:1} below).
A quick corollary, which is another exhibit of 
the similarities between $\omega_1$ and Suslin trees, is the following.
Its proof is immediate since the road space of a tree satisfying the assumptions is locally metrizable and
a Suslin tree is $\omega_1$-compact (see Lemma \ref{lemgen:0} below).

\begin{cor}\label{cor:main1}
   Let $T$ be a tree of height $\omega_1$. 
   If $T$ contains an uncountable subset $S$ which is $\omega_1$-compact in the subspace topology
   (in particular, if $T$ is a Suslin tree),
   then its road space $R_T$ is non-contractible.
\end{cor}

Our second generalization of Deo and Gauld's theorem is to weaken the local metrizability of the target space $X$. 
First countability is not enough, as the cone $C\omega_1 = \omega_1\times[0,1] / (x,1) \sim (y,1)$ is first countable, contractible
and contains of course many copies of $\omega_1$ itself. 
(First countability follows from the easily proved fact that an open set in $\omega_1\times[0,1]$ containing $\omega_1\times\{1\}$ must contain a small strip
$\omega_1\times(a,1]$.)
What we were able to prove is the following.

\begin{thm}\label{thm:main3}
  Let $X$ be a regular space with $G_\delta$ points such that each point has a closed $[\aleph_0,\aleph_1]$-compact
  neighborhood.
  Let $S$ be a stationary subset of $\omega_1$ endowed with the subspace topology. 
  Then there is no continuous $h\colon S\times[0,1]\to X$ such that $h(\cdot,1)\colon S\to X$ has an uncountable image and $h(\cdot,0)$ is constant.
  In particular, if $X$ contains a copy of $S$ then $X$ is not contractible.
\end{thm}
Recall that a space is {\em $[\aleph_0,\aleph_1]$-compact} if and only if any
cover of $X$ of size $\le\aleph_1$ has a countable subcover. This is equivalent to
the property that given any uncountable subspace $E$ of $X$ there is $x\in X$ which is a {\em condensation point of $E$}, that is,
any neighborhood of $x$ contains uncountably many points of $E$.   
We do not know whether this theorem holds if $S$ is a Suslin tree but our main tool (any continuous map from $S$ to $X$
must be constant on the points above $x\in S$ if $x$ is high enough, see Lemma \ref{lemma:omega_1} below) does not, 
as we will show in another paper \cite{meszigueseventually}.
Notice that Theorem \ref{thm:main3} holds if $X$ is locally metrizable.
Indeed, a stationary subset of $\omega_1$ is $\omega_1$-compact,
$\omega_1$-compactness is preserved by continuous functions and
is equivalent to Lindel\"ofness in metrizable spaces; hence the image of $S\times[0,1]$ in $X$ is (locally) Lindel\"of.

We give two proofs of Theorem \ref{thm:main2}. One is entirely topological and takes a few pages, in great part because we take the time to first show
that many well known properties of Suslin trees also hold for $\omega_1$-compact subsets of trees of height $\omega_1$.
These generalizations are probably part of the folklore and no more than standard exercices, but we found
convenient to gather the proofs of most of them in Section \ref{section:trees}.
The property that is central to our proofs is that any continuous map from $S$ to a locally metrizable space
must be constant on the points above $x\in S$ if $x$ is high enough in the tree (see Lemma \ref{lemgen:1} (i) below).
This is a slight generalization of a result of Stepr\=ans \cite{StepransTrees}.
Our other proof is
by forcing, and is based on another argument of Stepr\=ans in the same paper. 
It is quite short but uses classical theorems on ccc forcing, some of the properties in Section \ref{section:trees} and (a weaker version of) 
Theorem \ref{thm:main3}. Since the proof of Theorem \ref{thm:main3} is very similar to the topological proof of Theorem \ref{thm:main2},
we decided to present our arguments in the following order: first the topological proof
of Theorem \ref{thm:main2}, then that of Theorem \ref{thm:main3} and finally the forcing argument for Theorem \ref{thm:main2}.
Both proofs of Theorem \ref{thm:main2} use heavily the tree structure.
We do not know whether the result can be generalized for arbitrary
non-metrizable $\omega_1$-compact spaces.

\begin{qu}\label{q:0}
  Is there a locally metrizable contractible space containing a non-metrizable $\omega_1$-compact subspace~?
  Is there a manifold with these properties~?
\end{qu}

We will show in another paper that if $X$ is countably compact, non-compact and Type I, then $X$ is not contractible. 
A consequence is that Quesion \ref{q:0} has a negative answer for 
$\omega_1$-compact locally compact spaces under the proper forcing axiom  
{\bf PFA}. 

Looking only at road spaces of trees, we do not know whether contractibility entails $\R$-specialness 
(but actually have a strange feeling of having overlooked 
something simple):

\begin{q}\label{q:1}
  Let $T$ be a tree and $R_T$ be its road space.\\
  (a) Does the contractibility of $R_T$ implies that $T$ is $\R$-special~?\\
  (b) Does the non-contractibility of $R_T$ imply that $T$ contains a Suslin subtree~?
\end{q}

Notice that it is not possible for (a) and (b) to have both a positive answer since
there are models of set theory with a non-$\R$-special tree $T$ 
that does not contain a Suslin subtree, see \cite{Schlindwein:2003}. We do not know whether these trees have contractible road spaces.

We end this introduction with definitions and notations.
Any function in this note is assumed to be continuous otherwise stated.
A {\em contraction} is a function $h\colon X\times[0,1]\to X$ such that $h(x,1)=x$ and $h(x,0)$ is a constant map.
A space is {\em contractible} if there is a contraction $h\colon X\times[0,1]\to X$.
We often write $h_t(x)$ for $h(x,t)$. 

Recall that
a {\em tree} $T$ is a partially ordered set such that each point has a well ordered set of predecessors.
We define the {\em height} of $x\in T$ and of $T$, the $\alpha$-th-level $\lev_\alpha(T)$, the {\em chains} and {\em antichains} in $T$ as 
usual, see for instance \cite[Section II.5]{kunen} or \cite{Nyikos:trees}. A subset $E$ of $T$ is {\em order-dense}
iff for any $y\in T$ there is $y\in E$ with $x<y$.
Recall that all trees are endowed with the order (also called interval) topology. 
We often abbreviate `closed and unbounded' by {\em club}.
A tree is {\em rooted} iff it has a unique minimal element called the {\em root}.
A {\em subtree} $S$ is a subset of $T$ with order restricted to $S$. Notice that the induced topology on a subtree $S$ is finer (sometimes striclty)
than the one given by the 
order restricted to $S$. Both topologies agree if $S$ is closed in $T$.
We assume that our trees are Hausdorff, 
that is, if $x,y\in T$ are at a limit level and have the same predecessors, then $x=y$. (This could be false for a subtree.)
An {\em $\omega_1$-tree} has countable levels and height $\omega_1$.
A tree is {\em Suslin} if it has height $\omega_1$ and its chains and antichains are at most countable.
When $x\in T$ and $\alpha$ is an ordinal, write $T(x) = \{y\in T\,:\, y\ge x\}$,
$x\upharpoonright\alpha$ for the unique predecessor of $x$ at level $\alpha$ (if $x$ is below the $\alpha$-th level, $x\upharpoonright\alpha=x$)
and $T_{\le\alpha}$ for the subset of elements at level $\le\alpha$. 
If $E\subset T$, set $E^\downarrow =\{x\in T\,:\,\exists y\in E\text{ with }x\le y\}$ to be its downward closure.
We say that the tree $T$ is {\em $\R$-special} iff there is a strictly increasing (not necessarily continuous)
function $T\to \R$. Recall that $\R$-special $\omega_1$-trees exist in {\bf ZFC}.

The road space $R_T$ of a tree $T$ is obtained by joining consecutive points by an interval $[0,1]$, 
with $0$ glued to the lowest point and $1$ to the highest. We extend the order in the obvious way.
When convenient we consider $T$ as a subset of $R_T$.
The topology on the interior of the added intervals is that of $(0,1)$.
For $x\in T\subset R_T$, in order for 
$R_T$ to be (locally) connected any open set containing $x$ must contain a small portion of each interval emanating from $x$. 
In order for the space to be locally metrizable (and hence first countable),
we take these portions uniformly as follows.
Denote by $_x[0,1]_y\subset R_T$
the interval between the two consecutive points $x,y\in T$. If $A\subset[0,1]$, then ${}_xA_y$ is understood as the 
corresponding subset of $_x[0,1]_y$.
For singletons we usually write $_x a_y$ instead of $_x\{a\}_y$.
If $x\in T$, denote by $s(x)$ the set of its immediate successors and set 
$W_{x,n} =\displaystyle \bigcup_{y\in s(x)} \,_x[0,1/n)_y$.
If $x\in s(z)$ with $z\in T$, a basis for the neighborhoods of $x$ is given by 
$\{ W_{x,n} \cup \,_z(1-1/n,1]_x\,:\,n\in\omega\}$.
A basic neighborhood of $x\in T$ at a limit level
is obtained by 
choosing some $z\in T$, $z<x$ and $n\in\omega$ and taking the segment $\{y\in R_T\,:\,z<y<x\}$ union each $W_{w,n}$ 
for those $w\in T$ with $z<w<x$.
This makes $R_T$ locally embeddable in $\R^2$ (as seen by induction).
The induced topology on $T\subset R_T$ is that of $T$, and $R_T$ is arc connected if and only if $T$ is rooted.

%%%%%%%%%%%%%%%%%%%%%%%%%%%%%%%%%%%%%%%%%%%%%%%

\section{A collection of facts on uncountable $\omega_1$-compact subsets of trees}\label{section:trees}

Notations: Given a tree $T$, if $A,B\subset T$ and $x\in T$, 
$A<B$ means that 
each member of $A$ is $<$ each member of $B$. We denote $\{x\} < A$ by $x<A$
and $T(x)\cap A$ by $A(x)$.
Notice that $\text{Lev}_\alpha(A^\downarrow) = \text{Lev}_\alpha(T) \cap A^\downarrow$.

\begin{lemma}\label{lemgen:0}
   Let $T$ be a tree of height $\omega_1$ and $S\subset T$ be uncountable and endowed with the subspace topology. Then the following hold.\\
   (a) The subspace topology on $S^\downarrow$ always agrees with the topology given by the induced order.
       If $S$ is closed in $T$, then the subspace topology agrees with the topology given by the induced order on $S$.
       \\
   (b) If $S$ is $\omega_1$-compact, then it intersects a stationary subset of levels of $S^\downarrow$ (and of $T$).\\   
   (c) An antichain is closed discrete in $T$, and a closed discrete subset of $T$ is an at most countable union of antichains.\\
   (d) If $S$ is $\omega_1$-compact then so is $S^\downarrow$.\\
   (e) If $S$ is $\omega_1$-compact and does not contain an uncountable chain, then $S^\downarrow$ is Suslin.\\
   (f) If $A\subset S\subset T$, then $A$ is a maximal antichain in $S$ if and only if it is a maximal antichain in $S^\downarrow$.
\end{lemma}
\proof
(a) and (f) are straightforward. For (b), if $S$ avoids a club set of levels of $S^\downarrow$ 
then taking one point between consecutive avoided levels (when available) yields an uncountable discrete subset which is closed in $S^\downarrow$ (and thus in $T$).
Item (c) is proved e.g. in \cite[Thm 4.11]{Nyikos:trees}, and (d)--(e) follow immediately from it.
\endproof

\begin{lemma}\label{lemgen:1}
   Let $T$ be a tree of height $\omega_1$ and $S\subset T$ be uncountable and $\omega_1$-compact in the subspace topology.
   Then the following hold.\\
   (a) $S^\downarrow$ is the disjoint union of a countable set, a Suslin tree and at most countably many copies of $\omega_1$.
       In particular, $S^\downarrow$ is an $\omega_1$-tree. \\
   (b) There is $\alpha\in\omega_1$ such that $|S(x)|=|S^\downarrow(x)| = \aleph_1$ when $x$ is above level $\alpha$.\\
   (c) If $C\subset S^\downarrow$ is uncountable, there is $x\in S$ such that $S^\downarrow(x)\subset C^\downarrow$ and $|S(x)|=|S^\downarrow(x)| = \aleph_1$.\\
   (d) $C^\downarrow$ is $\omega_1$-compact for any $C\subset S^\downarrow$.\\
   (e) If $C\subset S^\downarrow$ is club in $S^\downarrow$, then $\{\gamma\,:\,\text{Lev}_\gamma(C^\downarrow)\subset C\}$ is club in $\omega_1$.\\
   (f) If $C\subset S^\downarrow$ is club in $S^\downarrow$, 
       then $S\cap C$ intersects a stationary subset of levels.\\
   (g) If $E_n\subset S^\downarrow$ are closed and order-dense in $S^\downarrow$ 
       for $n\in\omega$, then $\cap_{n\in\omega}E_n$ is also closed and order-dense.
       \\
   (h) Let $E,F\subset S\subset T$ be closed in $S$.
       If $|E\cap F| \le\aleph_0$, then $|E^\downarrow\cap F^\downarrow|\le\aleph_0$.\\
   (i) Let $f:S\to Y$ be continuous where $Y$ is a metrizable space. Then there is $\beta\in\omega_1$
       such that $f(S(x))$ is a singleton whenever $x$ is above level $\beta$ in $T$.
       In particular, $f$ has a countable image.
\end{lemma}
\proof
We will use several times (without acknowledging it explicitly) the fact that an at most countable intersection of club
subsets of $\omega_1$ is club.

(a) If $S$ contains an uncountable branch $B$, then there is some minimal $x_B\in B$ such that $S(x_B)$ is linearly ordered.
       Indeed, otherwise $S$ contains an uncountable antichain (take points branching away from $B$ above each height)
       and hence an uncountable closed discrete subset by Lemma \ref{lemgen:0} (c).
       Since the minimal elements of $\{x_B\,:\,S(x_B)\subset S\text{ is an uncountable branch}\}$ is an antichain, 
       $S$ contains at most countably many disjoint branches.
       Notice that a maximal branch in $S^\downarrow$ contains an unbounded branch of $S$.
       Removing the branches above each $x_B$ in $S^\downarrow$, what remains is either countable or a Suslin tree by Lemma \ref{lemgen:0} (e).

(b) follows immediately from (a) and the equivalent statement for Suslin trees.

(c) If $C\cap B$ is uncountable for some branch $B\subset S^\downarrow$, we are over. 
     If not, by (a) and (b)
     we can assume that $S^\downarrow$ is a Suslin tree such that
     $|S(x)|=\aleph_1$ for each $x\in S^\downarrow$. 
     Then the result is well known (see, for instance, the claim in Theorem 2.1 in \cite{DevlinShelah}).

(d) is immediate by (a) since an uncountable (downward closed) subset of a Suslin tree is a Suslin tree.

(e) By (d), $C^\downarrow$ is $\omega_1$-compact. Fix $\alpha$ given by (b) and some $\beta>\alpha$.
    Since each uncountable maximal branch $B$ of $C^\downarrow$ is a copy of $\omega_1$, $C$ contains a club set 
    of levels of $B$. By (a) we may assume that $C^\downarrow$ is Suslin.
    By Lemma \ref{lemgen:0} (f), 
    for each $n\in\omega$ we may find a countable antichain $A_n\subset C$ which is maximal in $C^\downarrow$ and such that 
    $A_{n+1} > A_n$ and each member of $A_n$ is above height $\beta$. Let $\gamma = \sup_{n\in\omega}\sup\{\text{height}(x)\,:\,x\in A_n\}$.
    By construction $\text{Lev}_\gamma(C^\downarrow)$ is the set of limit points of $\cup_{n\in\omega}A_n$, hence
    since $C$ is closed $\text{Lev}_\gamma(C^\downarrow)\subset C$.
    This shows that $\{\gamma\,:\,\text{Lev}_\gamma(C^\downarrow)\subset C\}$ is unbounded in $\omega_1$,
    and closedness is obvious.

(f) If $C\cap B$ is unbounded for some maximal uncountable branch $B\subset S^\downarrow$, then 
       it is homeomorphic to a club subset of $\omega_1$.
       By Lemma \ref{lemgen:0} (b) $S\cap B \cap C$ is thus stationary. 
       If $C\cap B$ is bounded for each uncountable branch of $S^\downarrow$,
       we may assume by (a) that $S^\downarrow$ is a Suslin tree in which $C$ is unbounded.
       By (c) it follows that $S$ is unbounded in $C^\downarrow$ as well. 
       By (e) and Lemma \ref{lemgen:0} (b),
       $S$ intersects $C$ on a stationary set of levels.

(g) Closedness is immediate, and order-density follows immediately by (a) and the fact that 
    the result holds for Suslin trees and $\omega_1$.
\iffalse
       hence let us show that $\cap_{n\in\omega}E_n$ is order-dense.
       First, since each uncountable maximal branch $B$ of $S^\downarrow$ is homeomorphic to $\omega_1$,
       $\cap_{n\in\omega_1} E_n$ is club in $B$. By (a) we may now assume that $S^\downarrow$ is Suslin.
       Let $\alpha$ be given by Lemma \ref{lemgen:1} (a).
       Let $x\in S^\downarrow$. We may assume that $x$ is at level above $\alpha$.
       Fix $\sigma:\omega\to\omega$ such that $|\sigma^{-1}(\{n\})| = \aleph_0$ for each $n\in\omega$.
       Since each $E_n$ is order-dense, for each $n$ we may take an antichain $A_n\subset E_{\sigma(n)}$ which is maximal in $S^\downarrow$
       such that $x<A_n < A_{n+1}$. Since $S^\downarrow$ is Suslin there is some $y\in T$ which is the limit of a sequence $x_n>x$
       such that $x_n\in A_n$. It follows that $y\in\cap_{n\in\omega}E_n$.
\fi

(h) 
   Let $\wb{E},\wb{F}$ be the closures in $S^\downarrow$ of $E,F$.
   If $\wb{E}\cap\wb{F}$ is unbounded in $S^\downarrow$, by (f) $S\cap\wb{E}\cap\wb{F} = (S\cap\wb{E})\cap (S\cap\wb{F}) = E\cap F$
   is unbounded, a contradiction.
   Hence $\wb{E}\cap\wb{F}$ is bounded and thus disjoint above some level $\alpha$.
   It follows that $E$ and $F$ cannot be both unbounded in the same uncountable branch. 
   It is well known (see e.g. \cite[Thm 6.18]{Todorcevic:1984} or \cite[Thm 2.1]{DevlinShelah}) that if $A,B$ are disjoint closed sets in a Suslin tree,
   then $A^\downarrow\cap B^\downarrow$ is at most countable. 
   Together with (a), this shows that $|E^\downarrow\cap F^\downarrow|\le\aleph_0$.
   
   (i) Our proof is a slight adaptation of Stepr\=ans topological proof in \cite{StepransTrees} that a real valued map with domain a Suslin tree has a countable image. 
       Denote the distance in $Y$ by $\text{dist}(\cdot,\cdot)$.
       By (b) me may assume that $S(x)$ is uncountable for each $x$.
       Set $D(\epsilon) = \{ x\in S\,:\,\text{diam}(f (S(x))) \le \epsilon\}$, where diam stands for diameter, that is, the supremum
       of the distances between two points in a set.
       Assume for now that $D(\epsilon)$ is order-dense in $S$ (and thus in $S^\downarrow$) when $\epsilon >0$.
       Let $\wb{D(\epsilon)}$ denote the closure of $D(\epsilon)$ in $S^\downarrow$.
       By (g), $D = \cap_{n\in\omega,n>0}\wb{D(1/n)}$ is closed and order dense in $S^\downarrow$, hence
       by (f) $S\cap D$ intersects stationary many levels above each $x\in S$ (since $|S(x)|=\aleph_1$).
       Moreover, $D$ is upward closed in $S$. 
       Denote by $A$ the minimal elements of $D\cap S$. Then $A$ is an antichain of $S^\downarrow$,
       let $\beta$ be the supremum of the heights of its members.
       For each $x\in S$ above level $\beta$, the diameter of $f\bigl(S(x)\bigr)$ is $0$, hence  
       $f\bigl(S(x)\bigr)= f(\{x\})$ and the lemma is proved.\\
       To finish, we now prove that $D(\epsilon)$ is order-dense in $S$.
       Suppose that it is not the case and let $x\in S$ be such that 
       $\text{diam}\bigl(f(S(y))\bigr) >\epsilon$ for each $y>x$, $y\in S$.
       We build inductively antichains $A^\alpha_n$ ($n\in\omega,\alpha\in\omega_1$) such that the following hold.\\
         $\bullet$ $A_n^\alpha\subset S$ is maximal above $x$, that is, in $S^\downarrow(x)$,\\
         $\bullet$ $A_{n+1}^\alpha > A_n^\alpha > A_m^\beta$ for each $n,m\in\omega$ and $\alpha > \beta$,\\
         $\bullet$  If $u\in A_n^\alpha$, $v\in A_{n+1}^\alpha$, then $\text{dist}(f(u),f(v))\ge\epsilon/4$.\\
       Assume that $A_n^\alpha$ is defined.
       Set 
         $$E=\{ z\in S\,:\, z>A_n^\alpha\text{ and dist}(f(z),f(u))\ge\epsilon/4 \text{, where $u$ is the member of $A_n^\alpha$ below $z$}\}.$$ 
       It is enough to see that $E$ is order-dense in $S^\downarrow(x)$, since then we may put its minimal elements in $A_{n+1}^\alpha$.
       Let thus $w\in S^\downarrow(x)$, $w>u\in A_n^\alpha$.
       Up to going further up, we may assume that $w\in S$.
       If $\text{dist}(f(w),f(u))\ge\epsilon/4$, then $w\in E$.
       If not, then 
       $\text{dist}(f(w),f(u))<\epsilon/4$.
       Choose $v\in S(w)$ such that $\text{dist}(f(w),f(v))>\epsilon/2$ (which exists since we assumed $\text{diam}\bigl(f(S(y))\bigr) >\epsilon$ for each $y>x$).
       Then 
         $$\text{dist}(f(v),f(u)) \ge \text{dist}(f(v),f(w)) - \text{dist}(f(w),f(u)) >\epsilon /4,$$ 
       and $v\in E$.
       If $A_n^\gamma$ is chosen for each $n\in\omega$ and each $\gamma<\alpha$, set $A_0^\alpha$ to be an antichain
       in $S$, maximal in $S^\downarrow(x)$, whose members are all $>\cup_{n\in\omega,\gamma<\alpha}A_n^\gamma$. This defines $A_n^\alpha$
       for each $n\in\omega,\alpha\in\omega_1$.\\
       Set $\beta(\alpha)$ to be $\sup\{ \text{height}(y)\,:\,y\in\cup_{n\in\omega} A_n^\alpha\}$,
       let $C$ be the closure in $\omega_1$ of $\{\beta(\alpha)\,:\,\alpha\in\omega_1\}$,
       and $C'$ be its derived set (that is its limit points).
       By construction, if $y\in S,y>x$ and the height of $y$ in $S^\downarrow$ is in $C'$, then $f$ is not continuous at $y$
       as they is a sequence of points in $S$ converging to $y$ whose images are $\ge\epsilon/4$ apart.
       But by Lemma \ref{lemgen:0} (b) (and the fact that $S(x)$ is $\omega_1$-compact), there must be such an $y$, a contradiction. 
       This shows that $D(\epsilon)$ is order-dense and concludes the proof.
\endproof

\section{Proofs}

%%%%%%%%%%%%%%%%%%%%%%%

\subsection{Proof of Theorem \ref{thm:main1}}

\proof
Let $T$ be a rooted $\R$-special tree with root $r$ and let $f:T\to\R$ be a strictly increasing function. By replacing $f(x)$ by $\sup_{y<x}f(y)$ 
when $x$ is at limit levels, we may assume that $f$ is continuous. By composing with a strictly increasing function, we
may assume that the range of $f$ is contained in $[0,1]$.
We first start to define the contraction as a map $h:T\times[0,1]\to R_T$.
First, set $h(x,1) = x$ and 
then define $h$ such that $x$ travels downwards in $R_T$, starting to move exactly at time $t=f(x)$ and reaching
$y<x$ ($y\in T$) exactly at time $t=f(y)$. Since $f$ is strictly increasing, there is time available to cross the interval between
consecutive points.
In less readable formulas, let $x$ at level $\alpha$ be given.
If $f(x)\le t$, we set $h(x,t) = x$.
If $\beta<\alpha$ and $f(x\upharpoonright \beta)\le t\le f(x\upharpoonright \beta+1)$,
we set $h(x,t)$ to be ${}_{x\upharpoonright \beta} k _{x\upharpoonright \beta+1}$ where
   $\displaystyle k=
    \frac{t-f(x\upharpoonright \beta )}{f(x\upharpoonright \beta+1) - f(x\upharpoonright \beta)}$.
Finally, if $t\le f(r)$, we set $h(x,t) = r$.
It should be clear that $h$ is continuous (on $T\times[0,1]$) and that $h(x,0) = r$ for all $x\in T$.
Notice that $h$ has the property that at time $f(x)$, all of $T(x)$ is squished onto $x$.
Actually, each point on the tree starts to move exactly when all the points above it reach it, all at the same time.
This enables to extend easily the map to all of $R_T$, a point in the segment ${}_x[0,1]_{y}$ does not move until
$y$ reaches it, and then it follows it until the end.
This gives the required contraction.
\endproof

Notice that if one sets $j_t(x) = h(x,1-t)$ for $t\in[0,1]$, $j_t(x) = x$ when $t<0$
and $j_t(x) = r$ when $t>1$, then $j$ is actually a flow, that is $j_t\circ j_s = j_{t+s}$.

%%%%%%%%%%%%%%%%%%%%%%%%%%%%%%%%

\subsection{Topological proof of Theorem \ref{thm:main2}}

Our proof relies on simple consequences of the properties given in Lemma \ref{lemgen:1}, especially (i).
When $S$ is a subset of a tree $T$, the height of a point of $S$ is to be understood as its height in $T$.
Recall that $S(x) = T(x)\cap S$.

\begin{lemma}
   Let $X$ be a space containing a closed metrizable $G_\delta$ subset $B\subset X$, 
   $S$ be an uncountable $\omega_1$-compact subset of a tree of height $\omega_1$ and
   $f:S\to X$ be continuous.
   Then there is $\alpha\in\omega_1$ such that
   either $f(S(x))\cap B = \varnothing$ or $f(S(x))$ is a singleton whenever $x$ is at height $\ge\alpha$.
\end{lemma}
\proof
   By Lemma \ref{lemgen:1} (b), above some level each $S(x)$ is uncountable, we assume for simplicity that this holds for each $x\in S$.
   Let $U_n$ be open sets such that $\cap_{n\in\omega} U_n = B$.
   By Lemma \ref{lemgen:1} (h), there is $\beta\in\omega_1$ such that for each $n$ we have
   $$ \left(\,\left(f^{-1}(B))\right)^\downarrow\cap \left(f^{-1}(X-U_n)\right)^\downarrow \right) \, - T_{\le\beta} = \varnothing.$$ 
   
   It follows that if $x\in S$ is at level above $\beta$,
   either $f(S(x))\subset B$ or $f(S(x))\cap B = \varnothing$.
   Let $E = \{x\in S\,:\, \lev(x)\ge\beta\text{ and }f(S(x))\subset B\}$.
   Then $E$ is an upward closed subspace of $S$, in particular it is an $\omega_1$-compact subspace of $T$.
   By Lemma \ref{lemgen:1} (i), there is $\alpha\ge\beta$ such that 
   $f$ is constant on $S(x)$ whenever $x\in E$ is at level $\ge\alpha$.
   This proves the lemma.
\endproof

\begin{cor}
   \label{cor:steprlocmet}
   Let $X$ be a space and $U\subset X$ be open such that
   $\wb{U}$ is contained in a metrizable open $V\subset X$.
   Let $S$ be an uncountable $\omega_1$-compact subset of a tree of height $\omega_1$.
   Let $h:S\times[0,1]\to X$ be continuous.
   Then there is $\alpha\in\omega_1$ such that for each $t\in[0,1]$ and each $x\in S$ above level $\alpha$,
   either $h_t^{-1}(U) \cap S(x) = \varnothing$, or $h_t$ is constant on $S(x)$.
\end{cor}
Again, $h_t$ stands for $h(\cdot,t)$.
\proof
   Let $\{t_n\,:\,n\in\omega\}$ be a countable dense subset of $[0,1]$. Set $B=\wb{U}$.
   Since $B$ is closed in the metrizable subset $V$, it is a $G_\delta$.
   The previous lemma shows that there is some $\alpha$
   such that when $x$ is above level $\alpha$, either $h_{t_n}^{-1}(B)\cap S(x) = \varnothing$ or $h_{t_n}$ is constant on $S(x)$
   for each $n\in\omega$.
   The result follows by continuity.
\endproof

\proof[Proof of Theorem \ref{thm:main2}]
Let $h:S\times[0,1]\to X$ be continuous such that $h_0(x) = u_0\in X$ and $h_1:S\to X$ has uncountable image.
We recall that we can assume that $T\supset S$ has height $\omega_1$.
The set of $x\in S$ such that $S(x)$ has uncountable image under $h_1$ is uncountable and downward closed, hence by Lemma \ref{lemgen:1} (c) there is $x$
such that the image of $S(z)$ under $h_1$ is uncountable for each $z\ge x$.
Up to replacing $S^\downarrow$ by $S^\downarrow(x)$ we assume that this holds for all $z\in S^\downarrow$.
For $x\in S$, set 
$$
\tau(x)=\sup\{t\,:\,h_t\text{ is constant on }S(x)\}.
$$
Then $\tau$ is an increasing map $S\to\R$, however $\tau$ is a priori neither continuous nor strictly increasing.
We will show that there is a closed unbounded $C\subset\omega_1$ such that $\tau$ is strictly increasing on the subspace of members of $S$ at 
levels belonging to $C$. 
Such a subspace is uncountable (and thus unbounded) in $S$ by Lemma \ref{lemgen:1} (f).
Hence as a partially ordered space $S$ contains an $\R$-special tree and thus
(at least) an uncountable antichain (see for instance \cite[Thm 4.29]{Nyikos:trees}), a contradiction
with the fact that $S$ is $\omega_1$-compact.
\\
It is enough to show that for each $x\in S$, there is $\alpha$ such that $\tau(y)>\tau(x)$ whenever $y>x$ is at level $\ge\alpha$. Indeed,
since the levels of $S^\downarrow$ are countable by Lemma \ref{lemgen:1} (a) and $\tau$ is increasing, a simple induction provides $C$.
Let $x\in S$ be fixed.
By continuity, $h_{\tau(x)}$ is constant on $S(x)$ with value $u = h_{\tau(x)}(x)$ and thus $\tau(x)<1$ since $h_1(S(x))$ is uncountable. 
(While it is not needed, notice that $h$ restricted to the subspace $S(x) \times [\tau(x),1]$ 
'contracts' all of $S(x)$ to the point $u$.)
Since $X$ is locally metrizable, we may choose an open $U\ni u$ such that $B=\wb{U}$ is contained
in an open metrizable set.
Let $\alpha$ be given by 
Corollary \ref{cor:steprlocmet}. We may assume that $\alpha>\text{height}(x)$.
Assume that there is $y$ at level above $\alpha$ such that for each $t>\tau(x)$, $h_t$ is not constant on $S(y)$.
By definition of $\alpha$, this implies that $h_t(S(y))\cap U = \varnothing$ and
in particular that $h_t(y)\not\in U$ for each $t>\tau(x)$. But this contradicts the continuity of $h$ since $h_{\tau(x)}(y) = u\in U$.
Hence, $h_t$ is constant on $S(y)$ for at least one $t>\tau(x)$ and thus $\tau(y)>\tau(x)$.
This finishes the proof.
\endproof

%%%%%%%%

\subsection{Proof of Theorem \ref{thm:main3}}

The proof is almost exactly the same as that of Theorem \ref{thm:main2} once we have an equivalent of Corollary \ref{cor:steprlocmet} in this context.
This is given by Corollary \ref{cor:omega_1} below.
Lemma \ref{lemma:omega_1} below plays the role of Lemma \ref{lemgen:1} (i). 
We first state the following easy fact.

\begin{lemma}\label{lemma:staclub}
   Let $S\subset\omega_1$ be stationary, endowed with the subspace topology.
   Then an at most countable family of club subsets of $S$ has a club intersection. 
\end{lemma}
\proof
   A direct proof is quite easy, but notice that
   since
   $S$ is an $\omega_1$-compact subspace of the tree $\omega_1$, the result is also a consequence of Lemma \ref{lemgen:1} (f)--(g).
\endproof

\begin{lemma}\label{lemma:omega_1}
  Let $S$ be a stationary subset of $\omega_1$ endowed with the subspace topology.
  If $Y$ is regular, $[\aleph_0,\aleph_1]$-compact and has $G_\delta$ points, then
  any continuous $f: S \to Y$ is eventually constant, that is, there is $\alpha\in\omega_1$ such that $f(\beta) = f(\alpha)$ for each $\beta\ge\alpha$,
  $\beta\in S$.
\end{lemma}

\proof 
  We start by showing that there is some $c\in Y$ such that $f^{-1}(\{c\})$ is club in $S$.
  If $f(S)$ is countable, then this is immediate.
  We thus assume that $f(S)$ is uncountable.
  Since $Y$ is $[\aleph_0,\aleph_1]$-compact, $f(S)$ has a condensation point $c\in Y$. 
  Since $Y$ is regular and has $G_\delta$ points,  we may choose 
  open sets $U_n\ni c$, $n\in\omega$, such that $\cap_{n\in\omega}U_n = \cap_{n\in\omega}\wb{U_n} =\{c\}$.
  Since $c$ is a condensation point, $f^{-1}(\wb{U_n})$ is club in $S$ for each $n$, hence
  $f^{-1}(\{c\})=f^{-1}(\cap_{n\in\omega}\wb{U_n}) = \cap_{n\in\omega}f^{-1}(\wb{U_n})$ is club in $S$ by Lemma \ref{lemma:staclub}.
  \\
  Now, since $f^{-1}(Y - U_n)$ is closed, it must be bounded, otherwise it intersects $f^{-1}(\{c\})$.
  It follows that $f^{-1}(Y-\{c\}) = \cup_{n\in\omega} f^{-1}(Y - U_n)$ is bounded in $S$, say by $\alpha$.
  Hence $f$ is constant on $S$ above $\alpha$.  
\endproof

\begin{cor}\label{cor:omega_1}
   Let $S$ be a stationary subset of $\omega_1$ endowed with the subspace topology.
   Let $Y$ be a regular space with $G_\delta$ points. Let $U\subset Y$ be open such that
   $\wb{U}$ is $[\aleph_0,\aleph_1]$-compact. Let $h:S\times[0,1]\to Y$ be continuous.
   Then there is $\alpha\in\omega_1$ such that for each $t\in[0,1]$ 
   either $h_t^{-1}(U) \cap [\alpha,\omega_1)\cap S = \varnothing$, or $h_t$ is constant on $[\alpha,\omega_1)\cap S$.
\end{cor}
\proof
   Fix a countable dense subset $\{t_n\,:\,n\in\omega\}\subset[0,1]$.
   By Lemmas \ref{lemma:staclub}--\ref{lemma:omega_1} we may fix $\alpha$ such that for each $n$
   either $[\alpha,\omega_1)\cap S \cap h_{t_n}^{-1}(\wb{U})=\varnothing$ or $h_{t_n}$ is constant above $\alpha$.
   The result follows by continuity.
\endproof

The proof of Theorem \ref{thm:main3} can now be done exactly along the same lines as that of Theorem \ref{thm:main2}, we thus only give a sketch.
Set $\tau = \sup\{t\,:\,h_t\text{ is eventually constant}\}$. 
Then $h_\tau$ is eventually constant and $\tau < 1$ since $h_1$ has uncountable image. 
Fix $\alpha$ minimal in $S$ such that $h_\tau$ is constant above $\alpha$, take an open $U$ containing $h_\tau(\alpha)$
such that $\wb{U}$ is $[\aleph_0,\aleph_1]$-compact. By Corollary \ref{cor:omega_1}
this contradicts the continuity of $h$.

%%%%%%%%%%%%

\subsection{Proof of Theorem \ref{thm:main2} by forcing}

Chapter VII of \cite{kunen} is a convenient reference for this subsection. 
Let $\langle X,\tau\rangle$ be a topological space in the ground model $V$. Given a forcing extension $V[G]$ 
by a generic filter $G$, we denote by $\tau(G)$ the topology for $X$ in $V[G]$ with base $\tau$. Notice that
in general, $\tau\not=\tau(G)$ since new unions may appear in $V[G]$.
A function in the ground model $f:X\to Y$ which is continuous remains continuous in $V[G]$. 
Also, being a $1$-to-$1$ function is preserved.
Some properties of $X$ are preserved in any forcing extension, for instance metrizability and the separation axioms
$T_i$ $i\le 3\frac{1}{2}$ (see, e.g., \cite[Lemma 22]{BaumgartnerTall}).
We will force with a Suslin tree with reversed order. We gather in the next lemma the well known properties of such forcings we are going to need.
\begin{lemma}\label{lemma:forcingfacts}
   Let $T$ be a Suslin tree in the ground model $V$.
   Let  $V[G]$ be the extension by a generic filter $G$
   when forcing over $V$ with $T$ with the reverse order. Then the following hold.\\
   (a) Cardinals and cofinalities are preserved between $V$ and $V[G]$.\\
   (b) If $S\subset\omega_1$ is stationary in $V$, it remains so in $V[G]$.
   \\
   (c) If $f,B\in V[G]$ and $f:\omega\to B$, then $f\in V$ (``no new countable sets are added'').
\end{lemma}
As said above, a classical reference for the proofs is \cite{kunen}: (a) is Theorem VII.5.10, (b) is Exercices VII.(H1)--(H2)
and (c) is Theorem VII.8.4. 
\\
\proof[Proof of \ref{thm:main2}]
Again, we may assume that $T$ has height $\omega_1$.
Let $S\subset T$ and $h:S\times[0,1]\to X$ be as in the statement of the theorem, so $h_0$ is constant and $h_1$ 
has uncountable image. Pick one preimage in $S$ for each point in the image of $h_1$, this 
defines an uncountable subset $W$ of $S$. 
Up to replacing $S^\downarrow$ by $S^\downarrow(x)$ for some $x$, by 
Lemma \ref{lemgen:1} (c)--(d) we may assume that $W^\downarrow \supset S$ and $|S^\downarrow(y)|=\aleph_1$ for all $y\in S$.
By Lemma \ref{lemgen:1} (a) $S^\downarrow$ contains either a Suslin tree or 
an uncountable chain (or both). In the latter case Theorem \ref{thm:main3} shows that $h$ cannot be continuous, we
may thus assume that $S^\downarrow$ is a rooted Suslin tree.
We now force with $S^\downarrow$ with reversed order and let $G$ be a generic filter.
Then $\cup G$ is a copy of $\omega_1^{V[G]}$ which is equal to $\omega_1^V$ 
since no new countable sets are added by Lemma \ref{lemma:forcingfacts} (c).
Also, $\cup G$ meets every level of $S^\downarrow$ since $D_\alpha = \{x\in S\,:\, \text{height}(x)\ge\alpha\}$ is
order-dense for each $\alpha$ (and the enumeration is in the ground model). 
Moreover, $E_\alpha = \{x\in S\,:\, \exists y\in W\text{ with height}(y)\ge\alpha \text{ and } y\ge x\}$
is also order-dense, hence $(\cup G)\cap W$ is uncountable (in $V[G]$).
By Lemma \ref{lemgen:0} (b) $S$ intersects stationary many levels of $S^\downarrow$ in
$V$, hence it remains so in $V[G]$ by Lemma \ref{lemma:forcingfacts} (b).
$S\cap (\cup G)$ is thus homeomorphic to a stationary subset of $\omega_1$.
But since $X$ remains locally metrizable,
$h$ remains countinuous in $V[G]$ and $h(\cdot,0)$ remains constant in $V[G]$, by Theorem \ref{thm:main3}
$h(\cdot,1)$ cannot have an uncountable image.
This contradicts the fact that $(\cup G)\cap W$ is uncountable and $h(\cdot,1)$ remains $1$-to-$1$ on $W$.
\endproof

%%%%%%%%%%%%%%%%%%

\end{document}